\documentclass[graybox]{svmult}

\usepackage{amsmath,amssymb,psfig,epsfig,latexsym,here}

\usepackage{mathptmx}       % selects Times Roman as basic font
\usepackage{helvet}         % selects Helvetica as sans-serif font
\usepackage{courier}        % selects Courier as typewriter font
\usepackage{type1cm}        % activate if the above 3 fonts are
\usepackage{makeidx}         % allows index generation
\usepackage{graphicx}        % standard LaTeX graphics tool
\usepackage{multicol}        % used for the two-column index
\usepackage[bottom]{footmisc}% places footnotes at page bottom

\makeindex

\newcommand{\bsq}{{\vrule height .9ex width .8ex depth -.1ex }}
\newcommand{\eeq}{\end{equation}}
\newcommand{\beql}[1]{\begin{equation}\label{#1}}
\newcommand{\beq}{\begin{equation}}
\newcommand{\eqn}[1]{(\ref{#1})}
\newcommand{\ZZ}{{\mathbb Z}}

\newcommand{\Sstyle}{\textstyle}
\newcommand{\Tstyle}{\scriptstyle}
\newcommand{\US}{ \_\hspace{.2ex} }

\begin{document}

\title*{\Large\bf Descending Dungeons and Iterated Base-Changing}
\titlerunning{Descending Dungeons}

\author{David Applegate, Marc LeBrun and N. J. A. Sloane}

\institute{David Applegate \at
AT\&T Shannon Labs,
180 Park Ave., Florham Park, NJ 07932-0971,
\\
\email{david@research.att.com}
\and
Marc LeBrun \at
Fixpoint Inc.,
448 Ignacio Blvd. \#239, Novato, CA 94949,
\email{mlb@well.com}
\and
N. J. A. Sloane \at
AT\&T Shannon Labs,
180 Park Ave., Florham Park, NJ 07932-0971,
\\
\email{njas@research.att.com}
}

\maketitle
% November 8, 2006; revised February 4, 2007 and August 27, 2007

\textsc{To our friend and former colleague Peter Fishburn,
on the occasion of his 70th birthday.}

\vspace{\baselineskip}

\abstract{
For real numbers $a, b> 1$, let $a_{\Sstyle b}$
(also written as $a \US b$)
denote the result of interpreting $a$ in base $b$ instead of base $10$.
We define ``dungeons'' (as opposed to ``towers'') to be numbers of the form
$a \US b \US c \US d \US \ldots \US e$,
parenthesized either from the bottom upwards
(preferred) or from the top downwards.
Among other things, we show that the sequences of
dungeons with $n^{\rm th}$ terms 
$10 \_ 11 \_ 12 \_\ldots\US (n-1) \US n$
or
$n \US (n-1) \_\ldots\_12 \_ 11 \_ 10$
grow roughly like $10^{10^{n \log \log n}}$,
where the logarithms are to the base 10.
We also investigate the behavior
as $n$ increases of the sequence
$a \US a \US a \US\ldots\US a$, with $n$ $a$'s,
parenthesized from the bottom upwards.
This converges either to a single number
(e.g. to the golden ratio if $a= 1.1$),
to a two-term limit cycle (e.g. if $a = 1.05$) or else diverges
(e.g. if $a = \frac{100}{99}$).
}

\vspace{0.8\baselineskip}

\vspace{0.8\baselineskip}

\setlength{\baselineskip}{1.5\baselineskip}
\section{Introduction}\label{Sec1}
The starting point for this paper was the question:
what is the asymptotic behavior of the sequences
\begin{eqnarray}\label{EQN1}
{} & 10, ~ 10_{{\Sstyle 11}},~  10_{{\Sstyle 11}_{\Sstyle 12}}, ~ 10_{{\Sstyle 11}_{\Sstyle 12_{\Sstyle 13}}}, ~ \ldots \, ,
\nonumber \\
{} & 10, ~ {\Sstyle 11}_{\Sstyle 10},~  12_{{\Sstyle 11}_{\Sstyle 10}}, ~ 13_{{\Sstyle 12}_{\Sstyle 11_{\Sstyle 10}}}, ~ \ldots \, ,
\end{eqnarray}
where, for real numbers $a, b > 1$, $a_{\Sstyle b}$ (or,
more conveniently although less graphically, $a\US b$)
denotes the result of interpreting $a$ in base $b$ instead of base $10$?
That is, if $a$ is a real number $>1$, with decimal expansion
\beql{EQN2}
a = \sum_{i=- \infty}^{k} c_i 10^i , \quad
\mbox{~for~some~} k \ge 0,
\mbox{~all~} c_i \in \{0,1, \ldots, 9 \}, \mbox{~and~} c_k \ne 0\,,
\eeq
and $b$ is a real number $>1$,
then
\beql{EQN3}
a_{\Sstyle b} := a\US b := \sum_{i=- \infty}^k c_i b^i \,.
\eeq
We use text-sized subscripts in expressions like $a_{\Sstyle b}$ to
help distinguish them from symbols with ordinary subscripts.
The sum in \eqn{EQN3} converges, since 
\beql{EQN3a}
1 < a_{\Sstyle b} < 9b^{k+1} / (b-1) \,,
\eeq
and $a_{\Sstyle b}$ is well-defined if we agree to avoid decimal
expansions ending with infinitely many 9's.
This restriction is needed, since (for example)
$3_{\Sstyle b} =3$ for any $b > 1$, whereas
$$2.999 \ldots_{\Sstyle b} = 2 + \frac{9}{b}
+\frac{9}{b^2} + \frac{9}{b^3} + \cdots = 2 + \frac{9}{b-1} \neq 3
$$
unless $b=10$.
Equation \eqn{EQN3} is meaningful for some values of $a$ and $b \le 1$,
but to avoid exceptions we only consider $a, b > 1$. In this range
$a \US  b$ is a binary operation for which $10$ is
both a left and right unit.

In fact, since the iterated subscripts can be grouped either
from the bottom upwards or from the top downwards, 
there are really four sequences to be considered (it
is convenient to index these sequences starting at $10$):
\renewcommand{\arraystretch}{1.25}
$$
\begin{array}{lll}
(\alpha ) = (\alpha_{10}, \alpha_{11}, \alpha_{12}, \ldots) &:= & 10, \, 10\_11, \, 10\_(11\_12), \, 10\_(11\_(12\_13)), \, \ldots \,, \\ 
(\beta ) = (\beta_{10}, \beta_{11}, \beta_{12}, \ldots \,) &:= & 10, \, 10\_11, \, (10\_11)\_12, \, ((10\_11)\_12)\_13, \, \ldots \,, \\ 
(\gamma ) = (\gamma_{10}, \gamma_{11}, \gamma_{12}, \ldots \,) &:= & 10, \, 11\_10, \, 12\_(11\_10), \, 13\_(12\_(11\_10)), \, \ldots \,, \\ 
(\delta ) = (\delta_{10}, \delta_{11}, \delta_{12}, \ldots ~) &:= & 10, \, 11\_10, \, (12\_11)\_10, \, ((13\_12)\_11)\_10, \, \ldots \,.
\end{array}
$$
\renewcommand{\arraystretch}{1.0}
Sequence $(\alpha )$, for example, begins
$$
\begin{array}{l}
10, ~~10\_11 = 11, ~~10\_(11\_12) = 10\_13 =13 \,, \\ [+.15in]
10\_(11\_(12\_13)) = 10\_(11\_15) = 10\_16 = 16\,, \\ [+.15in]
%10, ~~10\_11 = 11, ~~10\_(11\_12) = 10\_13 =13, ~~
%10\_(11\_(12\_13)) = 10\_(11\_15) = 10\_16 = 16\,, \\ [+.15in]
10\_(11\_(12\_(13\_14))) = 10\_(11\_(12\_17)) = 10\_ (11\_19) = 10\_20 =20 , ~~ \ldots
\end{array}
$$
The terms grow quite rapidly---see Table \ref{T1}.
These are now sequences A121263, A121265, A121295 and A121296 in \cite{OEIS}.
\begin{table}[htb]
\caption{Initial terms of sequence $(\alpha )$, $(\beta )$, $(\gamma )$, $(\delta )$.}
\label{T1}
$$
\begin{array}{|r|r|r|r|r|} \hline
n & ( \alpha )~~~~~~~~~~~~~~~~~ & ( \beta )~~~~~~~~~ & ( \gamma )~~~~~~~~~~~~~~~~ & ( \delta )~~~~~~~~~~ \\ \hline
10 & 10 & 10 & 10 & 10 \\
11 & 11 & 11 & 11 & 11 \\
12 & 13 & 13 & 13 & 13 \\
13 & 16 & 16 & 16 & 16 \\
14 & 20 & 20 & 20 & 20 \\
15 & 25 & 30 & 25 & 28 \\
16 & 31 & 48 & 31 & 45 \\
17 & 38 & 76 & 38 & 73 \\
18 & 46 & 132 & 46 & 133 \\
19 & 55 & 420 & 55 & 348 \\
20 & 65 & 1640 & 110 & 4943 \\
21 & 87 & 11991 & 221 & 22779 \\
22 & 135 & 249459 & 444 & 537226 \\
23 & 239 & 14103793 & 891 & 11662285 \\
24 & 463 & 5358891675 & 1786 & 46524257772 \\
25 & 943 & 19563802363305 & 3577 & 1092759075796059 \\
\ldots & \ldots & \ldots & \ldots & \ldots \\
30 & 38959 & 3.6053\ldots \times 10^{80} & 171999 & 2.5841\ldots \times 10^{89} \\
\ldots & \ldots & \ldots & \ldots & \ldots \\
35 & 9153583 & 8.6168\ldots \times 10^{643} & 41795936 & 1.2327\ldots \times 10^{898} \\
\ldots & \ldots & \ldots & \ldots & \ldots \\
100 & 4.0033\ldots \times 10^{57} & \ldots & 4.9144\ldots \times 10^{114} & \ldots \\
\ldots & \ldots & \ldots & \ldots & \ldots \\
\ldots & 6.8365\ldots \times 10^{1098} & \ldots & 3.4024\ldots \times 10^{917} & \ldots \\
& ~\mbox{at}~ n = 109 & & ~\mbox{at}~ n = 103 & \\
%\ldots & 6.8365\ldots \times 10^{1098} ~\mbox{at}~ n = 109 & \ldots & 3.4024\ldots \times 10^{917} ~\mbox{at}~ n = 103 & \ldots \\
\hline
\end{array}
$$
\end{table}

In Theorem \ref{Th1} we will show that,
if $s_n$ is the $n^{\rm th}$ term in any of the four sequences $(\alpha)$,
$(\beta)$, $(\gamma)$ or $(\delta)$,
indexed by $n = 10, 11, \ldots$,
then 
\beql{EqNew1}
\log \log s_n ~\sim n ~\log \log n \mbox{~as~} n \rightarrow \infty
\eeq
(in this paper all logarithms are to the base 10).

Since expressions like
$$10^{{\Tstyle 11}^{\Tstyle{12}^{\Tstyle{13}}}}$$
are called {\em towers}, we will call expressions like those in
\eqn{EQN1} and $(\alpha)$, $(\beta)$, $(\gamma)$ or $(\delta)$, {\em dungeons}.
For reasons that will be given in \S\ref{Sec2}, we believe that the standard
parenthesizing of dungeons should be from the bottom upwards,
and we will take this as the default meaning if
the parentheses are omitted.
For towers of exponents, parenthesizing from the top downwards is clearly better
(for otherwise the tower collapses).
The tower with $n^{\rm th}$ term
$$
t_n :=
10 \uparrow (11 \uparrow (12 \uparrow \cdots ((n-1) \uparrow n)\cdots)) 
\,,~ n=10,11,\ldots \,,
$$
(where $a \uparrow b$ denotes $a^b$) has the property that
the iterated logarithm $\log^{(n)} t_n \to \infty$
(note that $\log^{(n)} t_n$ is well-defined for $n$ sufficiently large).
When parenthesized from the bottom upwards, the tower with $n^{\rm th}$ term
$$
u_n ~ := ~ (\cdots(( 10 \uparrow 11 ) \uparrow 12) \cdots (n-1))\uparrow n 
~ = ~ 10^{11 \cdot 12  \cdot \, \cdots \, \cdot n} \,,~ n=10,11,\ldots \,,
$$
has the property that $\log \log u_n \sim n \log n$.
Equation \eqn{EqNew1} shows that
the dungeon sequences have a slower growth rate
than either version of the tower.

In \S\ref{Sec3} and \S\ref{Sec4} we prove Theorem \ref{Th1} and give some other properties of these
sequences, such as the fact that sequence $( \alpha )$ 
converges $10$-adically---for example,
from a certain point on, the last ten digits are always
$\ldots9163204655$.

In \S\ref{Sec5} we investigate the behavior as $n$ increases of 
the sequence with $n^{\rm th}$ term ($n = 1, 2, \ldots$)
\beql{EqEE}
a(n) :=
a\US (a\US (a\US (a\US  \cdots a))) ~~~~
(\mbox{with $n$ copies of $a$}) \,
\eeq
for a fixed real number $a>1$.
If the parameter $a$ exceeds $10$ this sequence certainly diverges, and for
$a=10$ we have $a(n) = 10$ for all $n \ge 1$.
Somewhat surprisingly, it seems hard to
say precisely what happens for $1 < a < 10$.
The mapping from $a(n)$ to $a(n+1) = a_{\Sstyle {a(n)}}$
is a discrete dynamical system,
which converges either to a single number
(e.g. to the golden ratio if the parameter $a= 1.1$),
to a two-term limit cycle (e.g. if $a = 1.05$) or diverges
(e.g. if $a = \frac{100}{99}$).
But we do not have a simple characterization of 
the parameters $a$ that fall into the different classes.
 
Section \ref{Sec2} contains some general properties of the subscript notation.

The following definition will be used throughout.
If $a>1$ is a fixed real number with decimal expansion given by \eqn{EQN2}
and $x$ is any real number, we define the Laurent series
\beql{EQN5}
L^{\langle a \rangle}(x) := \sum_{i = - \infty }^{k} c_{i} x^{i},
\eeq
so that $a \US  b = L^{\langle a \rangle}(b) $.
We use angle brackets to show the dependence on
the parameter $a$.
Note also that $L^{\langle a \rangle}(10) = a_{{\Sstyle 10}} = a$ for all $a$.

%\noindent{\bf Remarks.}
\begin{remark}
The choice of base $10$ in this paper was a matter of
personal preference.
\end{remark}
\begin{remark}
To answer a question raised by some readers of an early
draft of this paper, as far as we know there is no
connection between this work and the
base-changing sequences studied by Goodstein \cite{Good}.
\end{remark}

\section{Properties of the subscript notation}\label{Sec2}
In this and the following section
we will be concerned with the numbers $a_{\Sstyle b}$ defined in \eqn{EQN3} when
$a$ and $b$ are integers $\ge 10$.

\begin{lemma}\label{LC}
Let $N= 
\sum_{i=0}^k \nu _i 10^i$, where the $\nu _i$ are nonnegative integers 
$($not necessarily in the range $0$ to $9)$, and suppose $b$ is an integer $\ge 10$.
Then
\beql{EqC1}
N_{\Sstyle b} \ge 
\sum_{i=0}^k \nu _i b^i \,.
\eeq
\end{lemma}

\vspace*{+.1in}
\noindent{\bf Proof.}
If the $\nu _i$ are all in the range $\{0, \ldots , 9\}$
then the two sides of
\eqn{EqC1} are equal.
Any $\nu _i \ge 10$, say $\nu _i = 10 q + r$, $q \ge 1$,
$r \in \{0, \ldots , 9 \}$,
causes the term $\nu _i b^i$ on the right-hand side of \eqn{EqC1}
 to be replaced by $qb^{i+1} + rb^i \ge (10q +r) b^i = \nu _i b^i$
on the left-hand side, and so the difference
between the two sides can only increase.~~~$\bsq$

\begin{corollary}\label{C1}
If $f(x)$ is a polynomial with nonnegative integer coefficients,
and $b$ is an integer $\ge 10$, then $f(10)_{\Sstyle b} \ge f(b)$.
\end{corollary}

\begin{lemma}\label{L1}
Assume $a,b,a' , b'$ are integers $\ge 10$.
Then
\begin{itemize}
\item[$($i$)$]
$~~a' \ge a \mbox{~if~and~only~if~} a'_{\Sstyle b} \ge a_{\Sstyle b}$ \,,
\item[$($ii$)$]
$~~b' \ge b \mbox{~if~and~only~if~} a_{{\Sstyle b}'} \ge a_{\Sstyle b}$ \,,
\item[$($iii$)$]
$~~(a+a')_{\Sstyle b} \ge a_{\Sstyle b} + a'_{\Sstyle b}$ \,,
\item[$($iv$)$]
%$~~a_{({\Sstyle b}+{\Sstyle b}')} \ge a_{\Sstyle b} + a_{{\Sstyle b}'}$ \,,
$~~a_{{\Sstyle (b+b')}} \ge a_{\Sstyle b} + a_{{\Sstyle b}'}$ \,,
\item[$($v$)$]
$~~a_{\Sstyle b} \ge \max \{a,b\}$ \,.
\end{itemize}
\end{lemma}

\vspace*{+.1in}
\noindent{\bf Proof.}
(i)~Suppose $a' = \sum_{i=0}^{r'} c'_i 10^i > a = 
\sum_{i=0}^{r} c _i 10^i$, with all $c'_i$, $c_i \in \{0, \ldots , 9\}$,
and let $k$ be the largest $i$ such that $c'_i \neq c_i$.
Then $a'_{\Sstyle b} - a_{\Sstyle b} =
\sum_{i=0}^k (c'_i -c_i ) b^i \ge b^k - \sum_{i=0}^{k-1} 9b^i > 0$.
The converse has a similar proof.
Claims (ii), (iv) and (v) are immediate, and (iii) follows from Lemma \ref{LC}.~~~$\bsq$

\vspace*{+.1in}
Note that all parts of Lemma \ref{L1}
may fail if we allow $a$ and $b$ to be less than 10
(e.g. $12_{\Sstyle 2} = 4 < 7_{\Sstyle 2} =7$;
$6_{\Sstyle 3} = 6 \ge 6_{\Sstyle 4} = 6$, but $3<4$).

\begin{lemma}\label{L2}
Assume $a,b,c$ are integers $\ge 10$.
Then
\beql{EqC2}
(a\US b)\US c \ge a\US (b\US c) \,.
\eeq
\end{lemma}

\vspace*{+.1in}
\noindent{\bf Proof.}
The left-hand side of \eqn{EqC2} is (in the notation of \eqn{EQN5})
$L^{\langle a \rangle} ( L^{\langle b \rangle}(10))_{\Sstyle c} =
(L^{\langle a \rangle} \circ L^{\langle b \rangle} ) (10)_{\Sstyle c}$,
where $\circ$ denotes composition.
The right-hand side is 
$L^{\langle a \rangle} ( L^{\langle b \rangle}(c)) =
(L^{\langle a \rangle} \circ L^{\langle b \rangle}) (c)$,
and the result now follows from
Corollary~\ref{C1}.~~~$\bsq$

\vspace*{+.1in}
We can now explain why we prefer the ``bottom-up'' parenthesizing of dungeons.
The reason can be stated in two essentially equivalent ways.
First, $a\US  (b\US (c\US d))$, say, is simply
$$L^{\langle a \rangle} \circ L^{\langle b \rangle} \circ L^{\langle c \rangle} (d) \, ,$$
whereas no such simple expression holds for $((a\US b) \US c)\US d$.
To put this another way, consider evaluating the $n^{\rm th}$
term of sequence $(\alpha )$ of \S\ref{Sec1}.
To do this, we must repeatedly calculate values of 
$r_{\Sstyle s}$ where $r$ is $\le n$ and $s$ is huge.
But to find the $n^{\rm th}$ term of $(\beta )$,
we must repeatedly calculate values of $r_{\Sstyle s}$ where $r$ is huge and $s \le n$.
The latter is a more difficult task,
since it requires finding the decimal expansion of $r$.
Again, when computing the sequence $a(1), a(2), a(3), \ldots$
for a given values of $a$ (see \eqn{EqEE}),
as long as the terms are parenthesized from the bottom
upwards, only one decimal expansion (of $a$ itself) is ever needed.

In \S\ref{Sec3} we will also need numerical estimates of $a_{\Sstyle b}$.
If $a,b \ge 10$ then $a_{\Sstyle b}$ is roughly $10^{\log a \log b}$ (remember that all
logarithms are to the base 10).
More precisely, we have:
\begin{lemma}\label{L3}
Assume $a, b$ are integers $\ge 10$.
Then
\beql{EqC3}
10^{\lfloor \log a \rfloor \lfloor \log b
\rfloor} \le 10^{\lfloor \log a \rfloor \log b} \le
a_{\Sstyle b} \le 10^{\log a \log b} \, .
\eeq
\end{lemma}

\vspace*{+.1in}
\noindent{\bf Proof.}
Suppose
$a = \sum_{i=0}^{k} c_i 10^i$
where
$k := \lfloor \log a \rfloor,
c_i \in \{0, 1, \ldots, 9\}$ for $i=0,1, \ldots, k$,
$c_k \ne 0$.
The left-hand inequalities in \eqn{EqC3} are immediate.
For the right-hand inequality we must show that
$$
\sum_{i=0}^{k} c_i b^i \le b^{\log a} \,,
$$
or equivalently that
$$
\log \{
c_k b^k (1+
\sum_{i=0}^{k-1}
\frac{c_i}{c_k b^{k-i}}
)
\}
\le
(\log b) \, \Big(
\log \{ c_k 10^k (1 + \sum_{i=0}^{k-1} \frac{c_i}{c_k 10^{k-i}}) \}
\Big) \,,
$$
and this
is easily checked to be true using $b \ge 10$.~~~$\bsq$

\section{Growth rate of the sequences 
$(\alpha )$, $(\beta )$, $( \gamma )$, $(\delta )$}\label{Sec3}

\begin{theorem}\label{Th1}
If $s_n$ $(n \ge 10 )$ denotes the $n^{\rm th}$ term in any of the sequences
$(\alpha )$, $(\beta )$, $( \gamma )$, $(\delta )$ then
$$\log \log s_n \sim n \log \log n ~~\mbox{as}~~ n \to \infty \,.$$
\end{theorem}

\vspace*{+.1in}
\noindent{\bf Proof.}
From Lemma \ref{L3} it follows that
$$
\prod_{i=10}^{n} \lfloor \log i \rfloor \le \log s_n \le \prod_{i=10}^{n} \log i \,.
$$
For the upper bound, we have
$$
\log \log s_n \le \sum_{i=10}^{n} \log \log i  \le n \log \log n \,.
$$
For the lower bound, 
$$
\log s_n \ge \prod_{i=10}^{n} \lfloor \log i \rfloor \ge
\prod_{i=10}^{n} \log i (1-\frac{1}{\log i}) \,,
$$
$$
\log \log s_n \ge \sum_{i=10}^{n} \log \log i  
- \sum_{i=10}^{n} \frac{1}{\log i} \,,
$$
and the right-hand side is
$\sim n \log \log n + O(n)$.~~~$\bsq$

A slight tightening of this argument
shows that there are positive constants $c_1, c_2$ such that
$$
n \log \log n - c_1 \frac{n}{ \log n} < \log \log s_n 
< n \log \log n - c_2 \frac{n}{ \log n}
$$
for all sufficiently large $n$.

\vspace*{+.1in}
Table \ref{T1} suggests that sequences
$(\beta )$ and $(\delta )$ grow faster than $(\alpha )$ and $(\gamma )$.
We can prove three of these four relationships.

\begin{theorem}\label{Th2}
For $n \ge 10$, $\beta_n \ge \alpha_n$ and $\delta_n \ge \gamma_n$.
\end{theorem}

\vspace*{+.1in}
\noindent{\bf Proof.}
This follows by repeated application of Lemma \ref{L2}.~~~$\bsq$

\begin{lemma}\label{LBG}
If for some real number $k > 10$ we have
$a \ge k b$ and $\log c \ge \log k / (\log~k~-~1)$,
then $a \US c \ge k(c \US b)$.
\end{lemma}

\noindent{\bf Proof.}
From Lemma \ref{L3} and the assumed bounds, we have
\begin{eqnarray*}
  & a \US c & \ge 10^{ \lfloor \log a \rfloor  \log c}  \nonumber \\
  &      {} & \ge 10^{(\log a - 1) \log c} \nonumber \\
  &      {} & \ge 10^{(\log b + \log k - 1) \log c}  \nonumber \\
  &      {} & = 10^{(\log c)(\log k - 1)} 10^{\log b \log c} \nonumber \\
  &      {} & \ge k(c \US b) \, . \nonumber 
\end{eqnarray*}
~~~$\bsq$

\begin{theorem}\label{ThBG}
For $n \ge 10$, $\beta_n \geq \gamma_n$.
\end{theorem}

\noindent{\bf Proof.}
From Table \ref{T1}, this is true for $n \le 23$.  For $n > 23$, since
$\beta_{n+1} = (\beta_n) \US  {(n+1)}$ and
$\gamma_{n+1} = (n+1) \US {\gamma_n}$,
the previous lemma (with $k=10^4$) gives us the result by
induction.~~~$\bsq$

\section{$p$-Adic convergence of the sequence $(\alpha )$}\label{Sec4}

%\vspace*{+.1in}
For the next theorem we need a further lemma.
Let us say that a polynomial $f(x) \in \ZZ [x]$ is {\em $m$-stable},
for a positive integer $m$, if all its coefficients except
the constant term are divisible by $m$.
In particular, if $f(x)$ is $m$-stable, $f(x) \equiv f(0)$ (mod $m$).

\begin{lemma}\label{L4}
If the polynomial $f(x) \in \ZZ [x]$ is $m$-stable and the
polynomial $g(x) \in \ZZ [x]$ is $n$-stable,
then the polynomial $h(x) := f \circ g (x)$ is $mn$-stable.
\end{lemma}

\vspace*{+.1in}
\noindent{\bf Proof.}
If $f(x) := \sum_i f_i x^i$, $g(x) := \sum_j g_j x^j$, 
then
$h(x) = \sum_i f_i \left( \sum_j g_j x^j \right)^i = \sum_k h_k x^k$ (say).
When the expression for $h_k$ $(k > 0)$ is expanded as a sum of monomials,
each term contains both a factor $f_i$ for some $i > 0$ and a factor $g_j$ for some $j > 0$.~~~$\bsq$

\begin{theorem}\label{Th3}
The sequence $\alpha_{10}, \alpha_{11}, \alpha_{12}, \ldots$ converges $10$-adically.
\end{theorem}

\vspace*{+.1in}
\noindent{\bf Proof.}
We know from the above discussions that, for any $10 \leq k < n$, 
$$
\alpha_n = {\Phi}^{[k]} ((k+1)\US (k+2)\US (k+3)\US \ldots\US n)\,,
$$
where ${\Phi}^{[k]}(x)$ is the polynomial
$$
{\Phi}^{[k]}(x) :=
L^{\langle {10} \rangle} \circ L^{\langle {11} \rangle}
\circ L^{\langle {12} \rangle} \circ 
\cdots \circ L^{\langle k \rangle} (x) \,.
$$
(We would normally write ${\Phi}_k(x)$, but since there are
already two different kinds of subscripts in this paper,
we will use the temporary notation
${\Phi}^{[k]}(x)$ in this proof instead.)
Now $L^{\langle {20} \rangle}(x)$, $L^{\langle {21} \rangle}(x) , \ldots ,
L^{\langle {29} \rangle}(x)$ are 2-stable and $L^{\langle {50} \rangle}(x) ,
\ldots , L^{\langle {59} \rangle}(x)$ are 5-stable, so by Lemma \ref{L4},
${\Phi}^{[59]} (x)$ is
$10^{10}$-stable.  This means that for $n \ge 60$, $\alpha_n \equiv
{\Phi}^{[59]} (0)$ (mod $10^{10}$),
and so is a constant (in fact 5564023619) $\bmod ~10^{10}$.
Similarly, $L^{\langle {500} \rangle}(x)$, $L^{\langle {501} \rangle}(x) ,
\ldots , L^{\langle {509} \rangle}(x)$ are 5-stable, 
so $\alpha_n$ is a constant $\bmod~10^{20}$ for $n \ge 510$;
and so on.~~~$\bsq$

%\paragraph{Remark.}
\begin{remark}
The same proof shows that $\alpha_{10}, \alpha_{11}, \alpha_{12}, \ldots$
converges $l$-adically, for any $l$ all
of whose prime factors are less than $10$.
\end{remark}

\section{The limiting value of $a\US a\US a\US a\US  \ldots$}\label{Sec5}
In this section we consider the behavior of 
the sequence $a(1), a(2), a(3), \ldots$
(see \eqn{EqEE}) as $n$ increases,
for a fixed real number $a$ in the range $1 < a < 10$.
For example, we have the amusing identity
\beql{EQN9}
{\Sstyle 1.1}_{{\Sstyle 1.1}_{{\Sstyle 1.1}_{{\Sstyle 1.1}_{{\Sstyle 1.1}_{{\Sstyle 1.1}_{{\Sstyle 1.1}_{{\Sstyle .}_{{\Sstyle .}_{\Sstyle .}}}}}}}}}
~=~  \frac{1+\sqrt{5}}{2}  \,.
\eeq
The sequence \eqn{EqEE} is the
trajectory of the discrete dynamical system
$x \mapsto L^{\langle a \rangle}(x)$ when started at $x=a$.
(Since $L^{\langle a \rangle}(10) = a$, we could also start all trajectories at $10$.)

\begin{figure}[htbp]
\caption{Trajectory of $L^{<1.1>}(x)$ starting at $x=1.1$.}
\label{converge}
\begin{center}
\epsfig{file=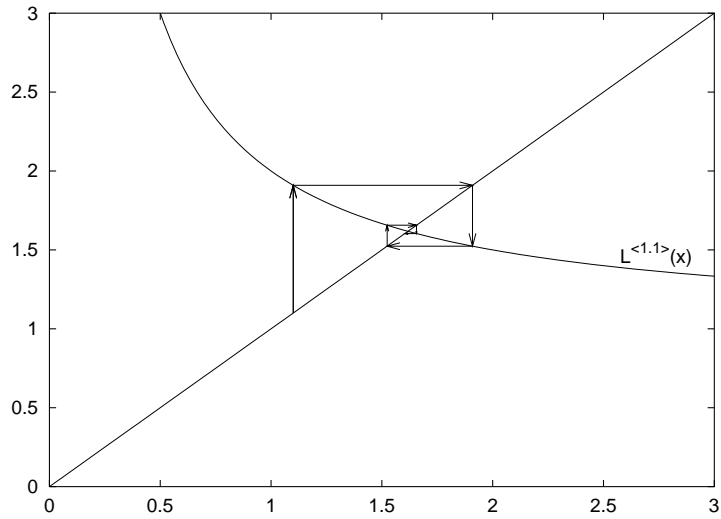,width=4in}
\end{center}
\end{figure}

Suppose $a = \sum_{i=0}^{ \infty} c_i 10^{-i}$ with 
all $c_i \in \{0,1, \ldots, 9 \}$ and $c_0 \ne 0$.
The graph of $y = L^{\langle a \rangle}(x)$
is a convex curve,
illustrated\footnote{This is a ``cobweb'' picture---compare Fig. 1.4 of \cite{Dev}.}
for $a=1.1$ in Figure~\ref{converge},
which decreases monotonically
from its value at $x=1$ (which may be infinite)
and approaches $c_0$ as $x \to \infty$. 
This curve therefore meets the line $y=x$ at a
unique point $x = \omega$ (say) in the range $x>1$.
The point $\omega$
is the unique fixed point for the dynamical system in the range of interest.

The general theory of dynamical systems \cite{Dev}, \cite{Lau} tells us that
the fixed point $\omega$
is respectively an attractor, a neutral point or a repelling point,
according to whether the value of the derivative 
${L^{\langle a \rangle}}'(\omega)$ 
is between $0$ and $-1$, equal to $-1$, or less than $-1$.
For our problem this does not tell the whole story, since we
are constrained to start at $a$.
However, since $L^{\langle a \rangle}(x)$ is a monotonically decreasing function,
there are only a few possibilities. Cycles of length three or more
cannot occur.

\begin{theorem}\label{Th4A}
For a fixed real number $a$ in the range $1<a<10$,
and an initial real starting value $x>1$, consider the trajectory
$x$, $L^{\langle a \rangle}(x)$,
$L^{\langle a \rangle} \circ L^{\langle a \rangle}(x)$,
$L^{\langle a \rangle} \circ L^{\langle a \rangle} \circ L^{\langle a \rangle}(x)$,
$\ldots$.
Then one of the following holds:
\begin{itemize}
\item[$($i$)$]
$~~x = \omega$ is the fixed point, and the
trajectory is simply $~\omega, \omega, \omega, \ldots$,
\item[$($ii$)$]
~~the trajectory converges to $\omega$,
\item[$($iii$)$]
$~~x$ is in a two-term cycle, and the
trajectory simply repeats that cycle,
\item[$($iv$)$]
~~the trajectory converges to a two-term limit cycle,
\item[$($v$)$]
~~the trajectory diverges,
alternately approaching $1$ and $\infty$.
\end{itemize}
\end{theorem}

\vspace*{+.1in}
\noindent{\bf Proof.}
If $a$ is an integer, then the trajectory is simply $x,a,a,a,\ldots$,
and either case (i) or (ii) holds.
Suppose then that $a$ is not an integer.
Since $a$ is fixed, we abbreviate
$L^{\langle a \rangle}$ by $L$ in this discussion, and write 
$L^{(k)}$ to indicate the $k$-fold composition of $L$, for $k=1,2,\ldots$.
Because $L(x)$ is strictly decreasing,
if $L^{(2)}(x) > x$, then $L^{(3)}(x) < L(x)$,
$L^{(4)}(x) > L^{(2)}(x) > x$;
if $L^{(2)}(x) < x$, then $L^{(3)}(x) > L(x)$,
$L^{(4)}(x) < L^{(2)}(x) < x$; and
if $L^{(2)}(x) = x$, then $L^{(3)}(x) = L(x)$,
$L^{(4)}(x) = L^{(2)}(x) = x$.
Hence if $x < L^{(2)}(x)$, then $x < L^{(2)}(x) < L^{(4)}(x) < \ldots$
and
if $x > L^{(2)}(x)$, then $x > L^{(2)}(x) > L^{(4)}(x) > \ldots$.
This means the even-indexed iterates form a monotonic sequence,
so either converge or are unbounded, and
similarly for the odd-indexed iterates.
Eq. \eqn{EQN3a} implies that if the trajectory diverges
then the lower limit must be $1$.~~~$\bsq$

Note also that
if $x < y < L^{(2)}(x)$, then $L^{(2k)}(x) < L^{(2k)}(y) < L^{(2k+2)}(x)$,
and if $x > y > L^{(2)}(x)$, then $L^{(2k)}(x) > L^{(2k)}(y) > L^{(2k+2)}(x)$.
So every $y$ between $x$ and $L^{(2k)}(x)$ converges to the same limiting
two-cycle as $x$ does, or diverges as $x$ does.

The following examples illustrate the five cases
in the situation which most interests us, the trajectory 
$a, a\US a, a\US (a\US a), \ldots$ of
\eqn{EqEE}, that is, when we set $x=a$ in the theorem.

(i) This case holds if and only if $a$ is one of $\{2, 3, \ldots, 9\}$.

(ii) Examples are
$a = 1+ \frac{m}{10}$, for $m \in \{1, \ldots, 9 \}$,
when $\omega = (1+ \sqrt{4m+1})/2$ is an attractor (see \eqn{EQN9});
$a = 1+ \frac{m}{100}$ for $m \in \{1, 2, 3 \}$,
when $\omega$, the real root $1.465\ldots$,
$1.695\ldots$ or $1.863\ldots$ of $x^3-x^2-m=0$
is an attractor;
and $a = 1+ \frac{4}{100}$, 
when $\omega = 2$ is neutral, but the trajectory still converges to $\omega$.

(iii) Examples are $a = 1+ \frac{m}{9}$, $m \in \{1, \ldots, 8 \}$,
 $\omega$ is a neutral point, and the two-term cycle is $\{ a, 10\}$.
(The trajectory does not include $\omega$.)

(iv) Examples are
$a = 1+ \frac{m}{100}$, $m \in \{5, \ldots, 9 \}$,
$\omega$ is a repelling point, and the trajectory approaches 
a two-term limit cycle consisting of a pair of solutions to
$L^{\langle a \rangle} \circ L^{\langle a \rangle}(x) = x$;
also $a = 1.1110000099$, $\omega$ is an attractor, but 
again the trajectory approaches a two-term cycle given by
$L^{\langle a \rangle} \circ L^{\langle a \rangle}(x) = x$.

(v) Examples are
$a = 1 + \frac{1}{10^r - 1}, r \in \{2,3,\ldots\}$,
$\omega$ is a repelling point, and the trajectory 
alternately approaches $1$ or $\infty$.

We do not know which values of $a$ fall into 
classes (ii) through (v). The distribution
of the five classes for $1 < a <10$ seems complicated.


\begin{thebibliography}{99}

\bibitem{Dev}
R. L. Devaney,
Dynamics of simple maps,
in {\em Chaos and Fractals}, ed. R.~L.~Devaney and L.~Keen,
Proceedings Symposia Applied Math., Vol. {\bfseries 39}, 
Amer. Math. Soc., Providence, RI, 1989.

\bibitem{Good}
R. L. Goodstein, On the restricted ordinal theorem,
{\em J. Symb.  Logic}, {\bfseries 9} (1944), 33--41.

\bibitem{Lau}
H. A. Lauwerier,
One-dimensional iterative maps,
in {\em Chaos}, ed. A. V. Holden, Princeton Univ. Press, 1986, pp. 39--57.

\bibitem{OEIS}
N. J. A. Sloane,
{\em The On-Line Encyclopedia of Integer Sequences},
published electronically at www.research.att.com/$\sim$njas/sequences/.

%\bibitem{SB}
%J. Stoer and R. Bulirsch,
%{\em Introduction to Numerical Analysis}, Springer-Verlag, NY, 1980.

\end{thebibliography}
\end{document}